\documentclass[11pt]{amsart}

\usepackage{amsxtra}
\usepackage{amssymb, amscd}

\newcommand\CC{{\mathbb C}}

\newcommand\RR{{\mathbb R}}
\newcommand\ZZ{{\mathbb Z}}

\newcommand\NN{{\mathbb N}}

\newcommand\SO{{\mathrm{SO}}}
\newcommand\SU{{\mathrm{SU}}}
\newcommand\U{{\mathrm{U}}}

\newcommand\End{\operatorname{End}}

\theoremstyle{plain}
\newtheorem{thm}{Theorem}[section]
\newtheorem*{thm*}{Theorem}

\newtheorem{lem}[thm]{Lemma}
\newtheorem{prop}[thm]{Proposition}
\newtheorem{cor}[thm]{Corollary}

\theoremstyle{definition}

\theoremstyle{remark}
\newtheorem{remark}[thm]{Remark}
\newtheorem*{remark*}{Remark}

\title[A sequence of matrix valued orthogonal polynomials]
{\sc A sequence of matrix valued orthogonal polynomials associated
to spherical functions}

\author{I. Pacharoni }
\author{P. Rom\'an}
\address{CIEM-FaMAF, Universidad Nacional de C\'ordoba, Ciudad
universitaria. CP. 5000 C\'ordoba. Argentina. Fax: 0054 351
4334054.} \email{pacharon@mate.uncor.edu}
\email{roman@mate.uncor.edu}

\thanks{This paper is partially supported by FONCyT grants PICT  03-10646, 03-13729,
 Secyt-UNC and by the ICTP Associate Scheme.}



\begin{document}

\begin{abstract} The main purpose of this paper is  to obtain an
explicit expression of a family of matrix valued orthogonal
polynomials $\{P_n\}_n$, with respect to a  weight $W$, that are
eigenfunctions of  a second order differential operator $D$. The
weight $W$ and the differential operator $D$   were found in
\cite{PT1}, using some aspects of the theory of the spherical
functions associated to the complex projective spaces.  We also find
other second order differential operator $E$ symmetric with respect
to  $W$ and we describe the algebra generated by $D$ and $E$.
\end{abstract}

\maketitle

\section{Introduction}

The theory of the harmonic analysis on homogeneous spaces is closely
connected with the theory of special functions. This is apparent,
for example, on the two dimensional sphere $S^2=\SO(3)/\SO(2)$,
where the harmonic analysis with respect to the action of the
orthogonal group is contained in the classical theory of the
spherical harmonics. In spherical coordinates the spherical
functions are the Legendre polynomials $P_n(\cos \theta)$. Also the
zonal spherical functions of the sphere $S^n=\SO(n+1)/\SO(n)$ are
given, in spherical coordinates, in terms of Jacobi polynomials
$P_n^{(\alpha,\alpha)}(\cos \theta)$, with $\alpha=(n-2)/2$. More
generally the zonal spherical functions on a Riemannian symmetric
space of rank one can always be expressed in terms of the classical
Gauss' hypergeometric functions, in the case of compact spaces we
get Jacobi polynomials.

As in the scalar case alluded above, in the matrix setting we also
have these three ingredients: the theory of matrix valued spherical
functions of any $K$-type, the matrix valued hypergeometric function
and the theory of matrix valued orthogonal polynomials. In this
paper we exhibit the interplay among these concepts in the case of
the complex projective space $P_n(\CC)=\SU(n+1)/\U(n)$.

The theory of matrix valued spherical functions goes back to
\cite{T0} and \cite{GV}, based on the foundational papers of
Godement and Harish-Chandra. In \cite{GPT1} we find  explicit
expressions for spherical functions, of any $K$-type associated to
complex projective plane $P_2(\CC)=\SU(3)/\U(2)$. This is
accomplished by associating to a spherical function $\Phi$ on $G$ a
 vector valued function $H$ defined on a complex affine plane
$\CC^2$, whose entries are given in terms of a special class of
generalized hypergeometric functions ${}_{p+1}\!F_p$.

The matrix valued hypergeometric function was studied  in \cite{T1}.
Let $V$ be a $d$-dimensional complex vector space, and let $A,B$ and
$C\in \End(V)$. The hypergeometric equation is
\begin{align}\label{hiper0}
z(1-z)F''(z) +(C-z (A+B+I))F'(z)- AB F(z)=0.
\end{align}
If the eigenvalues of $C$ are not in $-\NN_0$ we define the function
\begin{equation*}
{}_2\!F_1 \!\!\left( \begin{smallmatrix} A\,;\,B\\
C\end{smallmatrix} ; z\right)=\sum_{m=0}^\infty
\frac{z^m}{m!}(C;A;B)_m ,
\end{equation*}
where the symbol $(C;A;B)_m$ is defined inductively by
\begin{align*}
(C;A;B)_0 &=1, \\
(C;A;B)_{m+1}&=(C+m)^{-1}(A+m)(B+m)(C;A;B)_m, \quad m\geq 0.
\end{align*}
The function ${}_2\!F_1 \!\!\left( \begin{smallmatrix} A\,;\,B\\
C\end{smallmatrix} ; z\right)$ is analytic on $|z|<1$ with values in
$\End(V)$. Moreover if $F_0\in V$ then $F(z)= {}_2\!F_1 \!\!\left( \begin{smallmatrix} A\,;\,B\\
C\end{smallmatrix} ; z\right)\!F_0$ is a solution of the
hypergeometric equation \eqref{hiper0} such that $F(0)=F_0$.
Conversely any  solution $F$, analytic at $z=0$ is of this form.

The theory of matrix valued orthogonal polynomials, without any
consideration of differential equations goes back to \cite{K1} and
\cite{K2}. In \cite{D}, the study of the matrix valued orthogonal
polynomials which are eigenfunctions of certain second order
differential operators was started. The first explicit examples of
such  polynomials are given in \cite{GPT2},  \cite{GPT5} and
\cite{DG}.

Given a self adjoint positive definite matrix valued smooth weight
function $W=W(t)$ with finite moments, we can consider the skew
symmetric bilinear form defined for any pair of square matrix valued
polynomial functions $P(t)$ and $Q(t)$ by the numerical matrix
\[
( P,Q)= \int_{\RR}P(t) W(t)Q^*(t)dt,
\]
where $Q^*(t)$ denotes the conjugate transpose of $Q(t)$. This leads
to the existence of a sequence of matrix valued orthogonal
polynomials, that is a sequence $\{P_n(t)\}$, where $P_n$ is a
polynomial of degree $n$ with non singular leading coefficients and
$(P_n,P_m)=0$ if $n\neq m$.

We  also consider the skew symmetric bilinear form
\begin{equation}
\langle P,Q\rangle=(P^*,Q^*)^*,
\end{equation}

and we say that a differential operator $D$ is symmetric with
respect to $W$ if
\begin{equation}\label{sim}
\langle DP,Q\rangle=\langle P,DQ\rangle,
\end{equation}
for all matrix valued polynomial functions $P$ and $Q$.

Let $D$ be an ordinary linear differential operator with matrix
valued polynomial coefficients of degree less or equal to the order
of derivation.   If $D$ is symmetric with respect to $W$ then any
orthogonal sequence $\{P_n\}$, with respect to $(\cdot,\cdot)$,
satisfies
\begin{equation}\label{autofuncion}
DP_n^*=P_n^*\Lambda_n,
\end{equation}
for some numerical matrix  $\Lambda_n$.

Assume that the weight function $W=W(t)$ is supported in the
interval $(a,b)$ and let $D$ be a second order differential operator
of the form
\begin{equation}\label{D}
D=A_2(t)\frac{d^2}{dt^2}+A_1(t)\frac{d}{dt}+A_0(t),
\end{equation}
with matrix valued polynomial coefficients $A_j(t)$ of degree less
or equal to $j$.  In \cite{GPT5} (see also \cite{DG}) it is proved
that the condition of symmetry  for $D$ is equivalent to the
following three differential equations
\begin{equation}\label{partes}
\begin{split}
A_2^*W &= WA_2, \\
A_1^*W &= -WA_1 + 2(WA_2)', \\
A_0^*W &= WA_0 - (WA_1)' +(WA_2)'',\\
\end{split}
\end{equation}
with the boundary conditions
\begin{equation}\label{borde}
\lim_{t\to x}W(t)A_2(t)=0=\lim_{t\to
x}\big(W(t)A_1(t)-A_1^*(t)W(t)\big),\text { for $x=a,b$}.
\end{equation}

Finding explicit solutions of these equations is a highly non
trivial task. In \cite{DG} and \cite{DG2} the authors give some
families of examples. In \cite{PT1} one finds, for each dimension, a
three parameter family of pairs $\{W,D\}$ satisfying \eqref{partes}
and \eqref{borde}. These families arise from the representation
theory of Lie groups. After the change of variable $u=1-t$, the main
result in \cite{PT1} reads:

\begin{thm}\label{D1salto}  Let $\alpha,\beta >-1$, $0<k<\beta+1$ and  $\ell\in \NN$.
Let $ D$ be the differential operator defined by
$$ D=u(1-u)\frac{d^2}{du^2}+(C-uU)\frac{d}{du}-V,$$
with
\begin{align*}
C& =\sum_{i=0}^\ell (\beta+1+2i)E_{ii}+\sum_{i=1}^\ell iE_{i,i-1},
\quad
U=\sum_{i=0}^\ell (\alpha+\beta+\ell+i+2) E_{ii} , \displaybreak[0]\\
V&=  \sum_{i=0}^\ell i(\alpha+\beta+i-k+1)E_{ii}-
  \sum_{i=0}^{\ell-1} (\ell-i)(i+\beta-k+1)E_{i,i+1}.
\end{align*}
Then the differential operator $D$ is symmetric with respect to the
 weight matrix $W(u)=(1-u)^\alpha u^\beta Z(u)$  given by
$$ Z(u)=\sum_{i,j=0}^\ell\left( \sum_{r=0}^\ell
\textstyle \binom ri \binom rj \binom{\ell+k-r-1}{\ell-r}
\binom{\beta-k+r}{r} (1-u)^{\ell-r}u^{i+j}\right) E_{ij}. $$
\end{thm}

\begin{remark*}
Here, and in other parts of the paper, we use $E_{ij}$ to denote the
matrix with entry $(i,j)$ equal $1$ and $0$ otherwise.
\end{remark*}

This theorem is obtained from the first few steps in the explicit
determination of all matrix valued spherical functions associated to
the $n$-dimensional projective space $P_n(\CC)=\SU(n+1)/\U(n)$. The
idea, also used  in \cite{GPT1}, is to cook up from a matrix valued
spherical function a function $H$ which depends on a single variable
$u$. Using that the spherical functions are eigenfunctions of the
Casimir operator of $\SU(n+1)$ we deduce that, after an appropriate
conjugation, $H$ is an eigenfunction of an ordinary linear second
order matrix valued differential operator $D$. The fact that this
operator is symmetric with respect to the weight $W$ is a
consequence of the fact that the Casimir operator is symmetric with
respect to the $L^2$-inner product between matrix valued functions
on $\SU(n+1)$. At this point some readers may find useful to consult
references \cite{GV}, \cite{T0} and \cite{PT1}.

One of the main purposes  of this paper is to give explicit
expressions of  a sequence of orthogonal polynomials associated to
the weight  $W$ given in Theorem \ref{D1salto}. This is accomplished
by studying the vector space $V(\lambda)$ of all vector valued
polynomial solutions of the hypergeometric equation $DF-\lambda
F=0$. This space is non trivial if and only if
$$\lambda=\lambda_j(w)=-w(w+\alpha+\beta+\ell+j+1)-j(\alpha+\beta-k+1+j),$$
for some $w\in \NN_0$ and $j=0,1,\dots , \ell$. If the eigenvalues
$\lambda_j(w)$ are all different then there exists a unique
polynomial solution (up to scalars) of $DF=\lambda F$. In
Proposition \ref{dimVlambda} we compute, in the general case, the
dimension of the space $V(\lambda)$. With this knowledge at hand, we
construct a sequence of  polynomials $\{P_w\}$, by choosing the
$j$-th column of $P_w$ as a particular polynomial in
$V(\lambda_j(w))$. In Theorem \ref{orthopoly} we prove that
$\{P_w\}$ is an orthogonal sequence of matrix valued polynomials
such that $DP_w^*=P_w^*\Lambda_w(D)$, where $\Lambda_n(D)$ is the
real valued diagonal matrix
$$\Lambda_w(D)= \sum_{0\leq j\leq \ell} \lambda_j(w) E_{jj}.$$

The matrix  spherical functions associated to
$(G,K)=(\SU(n+1),\U(n))$ are eigenfunctions, not only of the Casimir
operator, but also of  any element in the algebra $D(G)^G$ of all
differential operators in $G$ which are left and right invariant
under multiplication by elements of $G$. In this case this algebra
is a polynomial algebra in $n$  algebraically independent
generators, one of them can be taken to be the Casimir operator of
$G$. For $n=2$, in \cite{GPT1} the explicit expression of this set
of generators was given and two differential operators $D$ and $E$
which commute were obtained. For a general $n$ we do not have simple
expressions for a complete set of generators of the algebra
$D(G)^G$, beyond the Casimir operator. However in this paper we are
able to find another second order differential operator $E$, which
commutes with $D$ and such that it is  symmetric with respect to $W$
(See Theorem \ref{E1salto}). The way in which we obtain this
operator is different to the one used in \cite{PT1} and it is
inspired in the operator $\tilde E$ given in \cite{RT}. Here we only
knew that such an operator should exist and after a trial and error
process we find it and prove that it is symmetric.

The sequence of matrix valued orthogonal polynomial constructed in
Theorem \ref{orthopoly} $\{P_w\}$ also satisfies  $EP_w^*=P_w^*
\Lambda_w(E)$ with
\begin{align*}
\Lambda_w(E)= \sum_{ j=0}^\ell &
(-w(w+\alpha+\beta+\ell+j+1)(\alpha-\ell+3j)\\
&\quad -j(j+\alpha+\beta-k+1)(\alpha+2\ell+3k)) E_{jj}.
\end{align*}

We also study the algebra generated by the differential operators
$D$ and $E$. In Theorem \ref{alggenDE} we prove that it is
isomorphic to the affine algebra of the following union of lines in
$\CC^2$:
$$\prod_{j=0}^\ell \left(y-(\alpha-\ell+3j)x+3j(\ell-j+k)(j+\alpha+\beta-k+1) \right).$$
Recently, in \cite{DIG} this situation is  considered  in the case
$\ell=2$. The authors conjecture that the algebra generated by $D$
and $E$ coincide with the algebra of all differential operators that
have the orthogonal polynomials $P_w$  as simultaneous
eigenfunctions.

\

\noindent {\bf Acknowledgement.}  We would like to thank Prof. Juan
Tirao for his continuous encouragement and for many useful comments
and suggestions that helped us to improve this paper.

\section{Orthogonal polynomials associated to the pair $\{W,D\}$}\label{orthogseccion}

 The aim of this section is to give explicitly a sequence of matrix valued
 orthogonal polynomials associated to the weight function $W$ and the differential operator $D$
introduced in Theorem \ref{D1salto}, i.e. we construct a sequence
$\{P_w\}$ of  orthogonal polynomials with respect to $W$, such that
$DP_w^*=P_w^*\Lambda_w$, where $\Lambda_w(D)$ is a real diagonal
matrix.

The columns $\{P_w^j\}_{j=0,\dots \ell}$ of $P_w^*$ are
$\CC^{\ell+1}$-valued polynomials such that $DP_w^j=\lambda_j(w)
P_w^j$ and $(P_w^j,P_{w'}^{j'})=\delta_{w,w'}\delta_{j,j'} n_{w,j}$,
for some positive real number $n_{w,j}$.

\subsection{Polynomial solutions of $DF=\lambda F$}
We start studying the $\CC^{\ell+1}$-vector valued polynomial
solutions of $DF=\lambda F$. We will find  all polynomials $F(u)$
such that
\begin{align}\label{hiperecuacion}
u(1-u)F''(u) +(C-u U)F'(u)- (V+\lambda) F(u)=0,
\end{align}
where the matrices $C,U,V$ are given in Theorem \ref{D1salto}. This
equation is an instance of a hypergeometric differential equation
studied in \cite{T1}. Since the eigenvalues of $C$ are not in
$-\NN_0$ the function $ F$ is characterized by $F_0=F(0)$. For
$\vert u\vert<1$ it is given by
\begin{equation}\label{hiper}
F(u)={}_2\!H_1 \!\!\left( \begin{smallmatrix} U\,;\,V+\lambda\\
C\end{smallmatrix} ; u\right)F_0=\sum_{i=0}^\infty
\frac{u^i}{i!}[C;U;V+\lambda]_i F_0,\qquad F_0\in \CC^\ell,
\end{equation}
where the symbol $[C;U;V+\lambda]_i$ is defined inductively by
\begin{align*}
[C;U;V+\lambda]_0 &=1, \\
[C;U;V+\lambda]_{i+1}&=(C+i)^{-1}\left(i(U+i-1)+V+\lambda\right)[C;U;V+\lambda]_i,
\end{align*} for all $i\geq 0$.

There exists  a polynomial solution of \eqref{hiperecuacion} if and
only if the coefficient $[C;U;V+\lambda]_i$ is singular for some
$i\in \ZZ$.
Let us assume that $[C;U;V+\lambda]_{w+1}$ is singular and that
$[C;U;V+\lambda]_w$ is not singular.

Since the matrix $(C+w)$ is invertible, we have that
$[C;U;V+\lambda]_{w+1}$ is singular if and only if
$(w(U+w-1)+V+\lambda)$ is singular. The matrix
$$M_w=(w(U+w-1)+V+\lambda)$$ is upper triangular and
$$(M_{w})_{j,j}=w(w+\alpha+\beta+\ell+j+1)+j(\alpha+\beta-k+1+j)+\lambda.$$
Therefore $[C;U;V+\lambda]_{w+1}$ is singular if and only if
\begin{equation}\label{autovalor}
\lambda=\lambda_j(w)=-w(w+\alpha+\beta+\ell+j+1)-j(\alpha+\beta-k+1+j),
\end{equation}
for some $0\leq j\leq \ell$.

We will distinguish the cases when the eigenvalues $\lambda_j(w)$
are all different (varying $j$ or $w$) or when they are repeated. We
start studying the polynomial solutions of \eqref{hiperecuacion} in
the first case.

\begin{prop} \label{unicoF0} Assume that all eigenvalues $\lambda_j(w)$ are
different.  If $\lambda=\lambda_j(w)$, for some $j=0,\dots, \ell$,
 then there exists a unique $F_0\in
  \CC^{\ell+1}$ (up to scalars) such that $F(u)={}_2\!H_1
  \left( \begin{smallmatrix} U\,;\,V+\lambda\\
C\end{smallmatrix} ; u\right)F_0$ is a polynomial function. Moreover
this polynomial is of degree $w$.
\end{prop}

\begin{proof} We have already observed that for
$\lambda=\lambda_j(w)=-w(w+\alpha+\beta+\ell+j+1)-j(\alpha+\beta-k+1+j)$,
the matrix $[C,U,V+\lambda]_{w+1}$ is singular. Then the function
$F(u)=\sum_{i=0}^\infty \frac{u^i}{i!}[C;U;V+\lambda]_i F_0$ is a
polynomial if and only if $F_0$ is a vector such that
\begin{equation}\label{F0condition}
[C,U,V+\lambda]_wF_0\in \ker (M_w);
\end{equation} where
$M_w=w(U+w-1)+V+\lambda_j(w)$.  The matrix $[C,U,V+\lambda]_w$ is
invertible, hence $F_0$ is univocally determined by an element in
the kernel of $M_w$. We have that
\begin{equation}\label{Mw}
M_w=\sum_{0\leq i\leq \ell}\bigl( (i-j)(\alpha+\beta-k+1+i+j+w)
E_{ii} - (\ell-i)(\beta-k+1+i)E_{i,i+1}\bigr).
\end{equation}
Since all eigenvalues $\lambda_j(w)$ are different we have that
$0\neq \lambda_j(w)-\lambda_i(w)=(i-j)(\alpha+\beta-k+1+i+j+w)$ if
$i\neq j$, hence the dimension of the kernel of $M_w$ is one.
Explicitly $(x_0,x_1, \dots, x_\ell)\in \ker(M_w)$ if and only if
\begin{equation}\label{kernelMw}
\begin{split}
 &x_i=\textstyle (-1)^{i+j}\binom{\ell-i}{\ell-j}
\frac{(\beta-k+1+i)_{j-i}}{(\alpha+\beta+j+i+w-k+1)_{j-i}} x_j\,
\qquad \text{ for }
i=0,\dots j, \\
&x_{j+1}= x_{j+2}=\cdots =x_{\ell}=0,
\end{split}
\end{equation}
where we use $(z)_r=z(z+1)\dots (z+r-1)$, $(z)_0=1$.

 Hence, up to
scalar, $F_0$ is uniquely determined by \eqref{F0condition} and it
is clear that $F(u)={}_2\!H_1 \!\!
\left( \begin{smallmatrix} U\,;\,V+\lambda\\
C\end{smallmatrix} ; u\right)F_0$ is a polynomial of degree $w$ with
 leading
coefficient  $\frac 1{w!}[C,U,V+\lambda_j(w)]_w F_0$. This completes
the proof of the proposition. \qed
\end{proof}

\smallskip
Now we have to study the case when some eigenvalues are repeated,
that is when there exist $w,w'\in \NN_0$ and $0\leq j,j'\leq \ell$
such that $\lambda_j(w)=\lambda_{j'}(w')$. We start observing the
following facts.

\begin{lem}\label{autovrepetido}
  If $\lambda_j(w)=\lambda_{j'}(w')$ for some $w,w'\in \NN_0$ and $0\leq j,j'\leq
  \ell$ then
  \begin{enumerate}
    \item [ i)] We have $w=w'$ if and only if $j=j'$.
    \item [ ii)] If $w'>w$ then $j>j'+1$.
  \end{enumerate}
\end{lem}
\begin{proof}
If $ \lambda_j(w)=\lambda_{j'}(w')$ then
\begin{align*}
 (w'-w)&(\alpha+\beta+\ell+1+w+w'+j')\\
 +& (j'-j)(\alpha+\beta-k+1+j+j'+w)=0.
\end{align*}

\noindent In particular if  $w'=w$, we have
$(j'-j)(\alpha+\beta-k+1+j+j'+w)=0$. We observe that $j\neq j'$
implies that $\alpha+\beta-k+1+j+j'+w>0$, because $\alpha>-1$,
$\beta-k+1>0$,
$j+j'\geq 1$ and $w\geq 0$.\\
Similarly if $j'=j$ we have $(w'-w)(\alpha+\beta+\ell+1+w+w'+j)=0$.
Since $\alpha>-1$, $\beta+\ell+1>0$ and  $w+w'+j\geq 1$ we obtain
that $(\alpha+\beta+\ell+1+w+w'+j)>0$ and therefore $w=w'$. This
completes the proof of i).

For ii) we start from
$$(w'-w)(\alpha+\beta+\ell+1+w+w'+j')= (j-j')(\alpha+\beta-k+1+j+j'+w),$$
and we observe that the left hand side of this identity, as well as
the factor $(\alpha+\beta-k+1+j+j'+w)$ are positive numbers, by
hypothesis, then we have $j>j'$. Finally suppose that $j=j'+1$ then
$(w'-w)(\alpha+\beta+\ell+w+w'+j)= (\alpha+\beta-k+w+2j),$
equivalently
$$(w'-w-1)(\alpha+\beta+\ell+w+w'+j)=-(w' + \ell-j+k).$$
The left  hand side is non negative while the right  hand side is
negative because $k>0$, which is a contradiction.
  \qed
\end{proof}

\

 Let $V(\lambda)$ be the  vector space  of all
$\CC^{\ell+1}$-vector valued  polynomials such that $DP=\lambda P$.
We observe that Proposition  \ref{unicoF0} said that if the
eigenvalues $\lambda=\lambda_j(w)$ are all different the dimension
of $V(\lambda)$ is one.
 The next proposition generalizes this result to the
 case when the eigenvalues $\lambda_j(w)$ are repeated.

\begin{prop}\label{dimVlambda} Let $\alpha,\beta>-1$, $0<k<\beta+1$ and let
$\lambda=\lambda_j(w)$, for some $w\in \NN_0$. Then
\begin{equation}\label{dimension}
\begin{split}
 \dim &\{P\in V(\lambda): \deg P\leq w\}\\
 &=\text{ card } \{ w': 0\leq w'\leq w\, , \,
\lambda=\lambda_{j'}(w'), \text{ for some } 0\leq j'\leq \ell\}.
\end{split}
\end{equation}
In particular
$$\dim V(\lambda)=\text{card }\{(w,j): \lambda=\lambda_j(w)\}.$$
\end{prop}
\begin{proof} We have already observed that for
$\lambda=\lambda_j(w)$ the function $F=F(u)$ is a polynomial
solution of $DF=\lambda F$ if and only if
$F(u)={}_2\!H_1(C,U,V+\lambda)F_0$ with $F_0\in \CC^{\ell+1}$ such
that $[C,U,V+\lambda]_wF_0\in \ker (M_{w,j}),$ where
$$M_{w,j}=\sum_{0\leq i\leq \ell}\bigl( (i-j)(\alpha+\beta-k+1+i+j+w)
E_{ii} - (\ell-i)(\beta-k+1+i)E_{i,i+1}\bigr)$$
We have that
$(i-j)(\alpha+\beta-k+1+i+j+w)\neq 0 $ if $i\neq j$. Hence  the
dimension of $\ker( M_{w,j})$ is one. Moreover it is generated by
$(x_0, \dots , x_\ell)\in \CC^{\ell+1}$ such that
\begin{equation}\label{kernelMw2}
\begin{split}
 &x_i=\textstyle (-1)^{i+j}\binom{\ell-i}{\ell-j}
\frac{(\beta-k+1+i)_{j-i}}{(\alpha+\beta+j+i+w-k+1)_{j-i}} \, \qquad
\text{ for }i=0,\dots j-1, \\
&x_j=1\\
 &x_{j+1}= x_{j+2}=\cdots =x_{\ell}=0,
\end{split}
\end{equation}
where we use $(z)_r=z(z+1)\dots (z+r-1)$, $(z)_0=1$.

If the eigenvalue $\lambda$ is repeated $s$ times and
$w_1=\min\{w\in \NN_0:\lambda=\lambda_j(w), 0\leq j\leq \ell\}$,
using Lemma \ref{autovrepetido}, we can assume that
$$\lambda=\lambda_{j_1}(w_1)=\cdots =\lambda_{j_s}(w_s)$$  with
$w_1<w_2<\cdots <w_s$ and  $j_1>j_{2}+1$, $j_2>j_{3}+1, \dots
,j_{s-1}>j_s$.

\smallskip
For $w=w_1$ and $j=j_1$ the matrix $[C,U,V+\lambda]_{w_1}$ is
invertible and $F_0$ is univocally determined by an element in $\ker
(M_{w_1,j_1})$, which is one dimensional, thus proving
\eqref{dimension} in this case.

Then to prove the proposition for any $w_r$ we proceed by induction
on $1\leq r\leq s$. Thus let us assume that for $2\leq r\leq s$ we
know that
$$\{P\in V(\lambda): \deg P\leq w_{r-1}\}=r-1.$$

Let $M_r=M_{w_r,j_r}$. As we remarked $0\neq P\in V(\lambda)$ is of
degree $w_r$ if and only if $P_0=P(0)$ satisfies $0\neq
[C,U,V+\lambda]_{w_r}P_0\in \ker (M_{r})$.

 Let
$$[C,U,V+\lambda]_{w_r}=N_{r}M_{{r-1}}\dots N_1M_{1}N_0,$$
where $N_i$ are invertible matrices. The leading coefficient $P_r$
of such a $P$ is uniquely determined, up to scalar, by the condition
$$M_rN_rM_{{r-1}}\dots N_1M_{1}N_0 P_0=0,$$ because we may assume that
$$P_r=N_rM_{{r-1}}\dots N_1M_{1}N_0 P_0=(x_0, \dots,
x_{j_r-1},1,0,\dots, 0).$$

Now let us prove that there exists $\tilde P\in V(\lambda)$ of
degree $w_r$, by constructing one by downward induction.

Let $v_r=(x_0, \dots, x_{j_r-1},1,0,\dots, 0)\in \ker(M_r)$ and let
$b_r=N_r^{-1}v_r$. The equation $b_r=M_{r-1}v_{r-1}$ has a unique
solution $v_{r-1}$ of the form $v_{r-1}=(z_0, \dots,
z_{j_r+1},0,\dots, 0)$ because $b_r=(y_0, \dots, y_{j_r+1},0,\dots,
0)$ with $y_{j_r+1}\neq 0$ and $M_{r-1}$ is upper triangular with a
unique zero in the main diagonal in the $j_{r-1}$-position.
Similarly let $b_{r-1}=N_{r-1}^{-1}v_{r-1}$, then there exists a
unique $v_{r-2}=(t_0, \dots, t_{j_r+2},0,\dots 0)$ such that
$M_{r-2}v_{r-2}=b_{r-1}$. In this way we construct the sequence
$v_r, v_{r-1}, \dots, v_0$ such that
\begin{align*}
 v_r &= N_rb_r= N_rM_{r-1}v_{r-1}= N_rM_{r-1}N_{r-1}M_{r-2}v_{r-2}=
\cdots \\
&= N_rM_{r-1}\dots N_1M_1N_0 v_0
\end{align*}
Hence $\tilde P={}_2\!H_1(C,U,V+\lambda)v_r$ is a polynomial in
$V(\lambda)$ of degree $w_r$.

Now we observe that
$$\{P\in V(\lambda): \deg P\leq w_r\}= \CC\tilde P \oplus
\{P\in V(\lambda): \deg P\leq w_{r-1}\}.$$ In fact it is clear that
the right hand side is a direct sum contained in the left hand side.
To prove the other inclusion we first observe that if $P\in
V(\lambda)$ and $\deg P<w_r$ then, as we saw, $\deg P\leq w_{r-1}$.
If $P\in V(\lambda)$ is of degree $w_r$ then the leading coefficient
of $P$ is equal to the leading coefficient of $t\tilde P$ for some
$t\in \CC$. Therefore $P-t\tilde P\in \{P\in V(\lambda): \deg P\leq
w_{r-1}\}$. This completes the proof of the proposition.\qed
\end{proof}

\subsection{Matrix valued orthogonal polynomials associated to
$\{W,D\}$.}\label{orthogsubseccion}

 We want to construct a
sequence $\{P_w\}_{w\geq 0}$ of  matrix valued orthogonal
polynomials with respect to the weight function $W$, with degree of
$P_w$ equal to $w$, with non singular leading coefficient  and that
satisfies $DP_w^*=P_w^*\Lambda_w$, where $\Lambda_w(D)$ is a real
diagonal matrix.

Then the columns $\{P_w^j\}_{j=0,\dots \ell}$ of $P_w^*$ are
$\CC^{\ell+1}$-valued polynomials such that $P_w^j$ and
$P_{w'}^{j'}$ are orthogonal to each other if $(j,w)\neq (j',w')$
and   they satisfy that $DP_w^j=\lambda_j(w) P_w^j$, where
$$\lambda_j(w)=-w(w+\alpha+\beta+\ell+j+1)-j(j+\alpha+\beta-k+1),$$
for  each $w\in\NN_0$, and $j=0,\dots,\ell$.

If an eigenvalue $\lambda=\lambda_j(w)$ is not repeated, then we
choose the unique $F_0\in \CC^{\ell+1}$ such that
\begin{equation}\label{F0}
[C,U,V+\lambda_j(w)]_w F_0=\sum_{0\leq i\leq j}\textstyle
(-1)^{i+j}\binom{\ell-i}{\ell-j}
\frac{(\beta-k+1+i)_{j-i}}{(\alpha+\beta+j+i+w-k+1)_{j-i}}\, e_i
\end{equation}
 where $e_i$ denotes
the $i$-th vector in the canonical basis of $\RR^{\ell+1}$. Then we
take
$$P_w^j(u)={}_2\!H_1 \left( \begin{smallmatrix} U\,;\,V+\lambda_j(w)\\
C\end{smallmatrix} ; u\right)F_0=\sum_{i=0}^\infty
\frac{u^i}{i!}[C;U;V+\lambda_j(w)]_i F_0$$ which is a polynomial
function of degree $w$ and satisfies $$DP_w^j(u)=\lambda_j(w)
P_w^j(u).$$ (See Proposition \ref{unicoF0}).

If an eigenvalue $\lambda=\lambda_j(w)$ is repeated we saw that,
$$\lambda= \lambda_{j_1}(w_1)= \lambda_{j_2}(w_2)=\cdots =
\lambda_{j_s}(w_s),$$  with $w_1<w_2<\cdots <w_s$ and $j_r\geq
j_{r+1}+1$, for $1\leq r\leq s-1$.

Let $V_r=\{ P\in V(\lambda): \deg P\leq w_r\}$, for $1\leq r\leq s$.
Then we saw, in the proof of Proposition \ref{dimVlambda}, that
$$0\neq V_1\subsetneq V_2\subsetneq \cdots \subsetneq V_s$$
with  $\dim V_s=s$ Now we take, for each $1\leq r\leq s$
$$0\neq P_{w_r}^{j_r}(u)={}_2\!H_1 \left( \begin{smallmatrix} U\,;\,V+\lambda_j(w)\\
C\end{smallmatrix} ; u\right)F_{0}^{j_r} \in V_r \text{  orthogonal
to }V_{r-1}.  $$

In this way, for each $w\in \NN_0$ we have defined $\ell+1$
orthogonal polynomial functions $P_w^0, P_w^1, \dots , P_w^\ell$ of
degree $w$.

\begin{thm}\label{orthopoly}
  Let $P_w(u)$ be the matrix whose rows are the vectors
$P_w^j(u)$. Then the sequence $\{P_w(u)\}_{w\in \NN_0}$ is an
orthogonal  sequence of matrix valued polynomials such that
$$DP_w^*(u)=P_w^*(u)\Lambda_w,$$ where  $\Lambda_w=\sum_{j=0}^\ell \lambda_j(w) E_{jj}$.
\end{thm}
\begin{proof}
Let $(w,j)\neq (w',j')$. If $\lambda_j(w)\neq \lambda_{j'}(w')$ then
$(P_w^j,P_{w'}^{j'})=0$ because $D$ is symmetric. If
$\lambda_j(w)=\lambda_{j'}(w')$ then $(P_w^j,P_{w'}^{j'})=0$ by
construction.  Therefore the matrices $P_w$ satisfies $(P_w,
P_{w'})=0$ if $w\neq w'$.

On the other hand we have that for each $w=0,1,2\dots $ the degree
of $P_w(u)$ is $w$ and the leading coefficient of $P_w$ is the non
singular triangular matrix
$$I+\sum_{ s< r}\textstyle (-1)^{r+s}\binom{\ell-s}{\ell-r}
\frac{(\beta-k+1+s)_{r-s}}{(\alpha+\beta+r+s+w-k+1)_{r-s}} E_{rs}.$$
This completes the proof of the theorem. \qed
\end{proof}

\section{The symmetry of the differential operator  $E$}\label{seccionE}

The aim of this section is to exhibit another second order ordinary
differential operator which is symmetric with respect to the weight
$W$.

\begin{thm} \label{E1salto} Let $\alpha, \beta>-1$, $0<k<\beta+1$ and $\ell\in \NN$.
Let $ E$ be the differential operator defined by
$$ E=(1-u)(Q_0+uQ_1)\frac{d^2}{du^2}+(P_0+uP_1)\frac{d}{du}-(\alpha+2\ell+3k)V,$$
with
\begin{align*}
Q_0&=\textstyle{\sum_{i=0}^\ell 3iE_{i,i-1}},\displaybreak[0]\\
Q_1&=\textstyle{\sum_{i=0}^\ell (\alpha-\ell+3i)E_{ii}},\displaybreak[0]\\
P_0&=\textstyle{\sum_{i=0}^\ell \big( (\alpha+2\ell)(\beta+1+2i)
-3k(\ell-i)-3i(\beta-k+i) \big)E_{ii}}\\
&\textstyle{\quad-\sum_{i=0}^\ell i(3i+3\beta-3k+3+\ell+2\alpha)E_{i,i-1}},\displaybreak[0]\\
P_1&=\textstyle{\sum_{i=0}^\ell
-(\alpha-\ell+3i)(\alpha+\beta+\ell+i+2)E_{ii}}\\
&\textstyle{\quad +\sum_{i=0}^\ell 3(\beta-k+1+i)(\ell-i)E_{i,i+1}},\displaybreak[0]\\
V&=  \textstyle\sum_{i=0}^\ell i(\alpha+\beta-k+1+i)E_{ii}-
  \sum_{i=0}^{\ell-1} (\ell-i)(\beta-k+1+i)E_{i,i+1}.
\end{align*}
Then $E$ is symmetric with respect to the
 weight matrix $W(u)=(1-u)^\alpha u^\beta Z(u)$, where $Z(u)$ is given  by
$$ Z(u)=\sum_{i,j=0}^\ell\left( \sum_{r=0}^\ell
\textstyle \binom ri \binom rj
\textstyle\binom{\ell+k-1-r}{\ell-r}\binom{\beta-k+r}{r}
(1-u)^{\ell-r}u^{i+j}\right) E_{ij}. $$
\end{thm}

\begin{proof}

We need to prove that the equations \eqref{partes} and \eqref{borde}
are satisfied. The equations in \eqref{partes} take the form
\begin{align}
   \label{eqI} &(Q_0^*+uQ_1^*)Z-Z(Q_0+uQ_1)=0,\displaybreak[0]\\
    \label{eqII}\begin{split} &(P_0^*+uP_1^*)Z+Z(P_0+uP_1)- 2Z( Q_1-Q_0-2uZQ_1)\\ &\qquad
    -    2(1-u)Z'(Q_0+uQ_1)-\tfrac{(\beta(1-u)-\alpha
    u)}{u}2Z(Q_0+uQ_1)=0,
    \end{split} \displaybreak[0]\\
    \label{eqIII} \begin{split}&
    P_1^*Z+(P_0^*+uP_1^*)Z'-Z'(P_0+uP_1)-ZP_1\\
    &\qquad +(\tfrac{\beta}{u}-\tfrac{\alpha}{1-u})\bigl(
    (P_0^*+uP_1^*)Z-Z(P_0+uP_1)\bigr)\\
    &\qquad -2(\alpha+2\ell+3k)(ZV-V^*Z)=0.
    \end{split}
\end{align}

The $ij$-entry in the left hand side of the equation \eqref{eqI} is
$$u(\alpha-\ell+3i)z_{ij}+3(i+1)z_{i+1,j} -
u(\alpha-\ell+3j)z_{ij}+3(j+1)z_{i,j+1}=0,$$
 because it is easy to verify  that
$$(i+1)z_{i+1,j}-(j+1)z_{i,j+1}=u(j-i)z_{ij}.$$

In order to prove the identity \eqref{eqII} we compute the
$ij$-entry of the matrices involved there:
\begin{align*}
&((P_0^*+uP_1^*)Z)_{ij}=u^{i+j}\sum_{r=\max(i,j)}^\ell\textstyle
\binom{r}{i}\binom{r}{j}\binom{\beta+r-1}{r}\binom{\ell+k-1-r}{\ell-r}(1-u)^{\ell-r}
\Big( (P_0)_{ii}\\ &\quad -
(r-i)(3i+3\beta+\ell+2\alpha-3k+6)-(\alpha-\ell+3i)(\alpha+\beta+\ell+i+2) \Big)\displaybreak[0]\\
& +u^{i+j}\sum_{r=\max(i-1,j-1)}^\ell \textstyle
\binom{r+1}{i}\binom{r+1}{j}\binom{\beta+r}{r+1}\binom{\ell+k-2-r}{\ell-r-1}(1-u)^{\ell-r}\\
&\,  \Big(
(r-i+1)(3i+3\beta+\ell+2\alpha-3k+6)+(\alpha-\ell+3i)(\alpha+\beta+\ell+i+2))\Big)\displaybreak[0]\\
& +u^{i+j}\!\!\!\!\sum_{r=\max(i-1,j)}^\ell \textstyle
\binom{r}{i-1}\binom{r}{j}\binom{\beta+r-1}{r}\binom{\ell+k-1-r}{\ell-r}(1-u)^{\ell-r}3(\ell-i+1)(\beta+i-k),
\end{align*}
\begin{align*}
&(Z(P_0+uP_1))_{ij}=u^{i+j}\sum_{r} \textstyle
\binom{r}{i}\binom{r}{j}\binom{\beta+r-1}{r}\binom{\ell+k-1-r}{\ell-r}(1-u)^{\ell-r} \Big( (P_0)_{jj}\\
&\quad -(r-j)(3j+3\beta+\ell+2\alpha-3k+6)-(\alpha-\ell+3j)(\alpha+\beta+\ell+j+2) \Big)\displaybreak[0]\\
&+ u^{i+j}\sum_{r} \textstyle
\binom{r+1}{i}\binom{r+1}{j}\binom{\beta+r}{r+1}\binom{\ell+k-2-r}{\ell-r-1}(1-u)^{\ell-r}\\
& \big( (r-j+1)(3j+3\beta+\ell+2\alpha-3k+6)+(\alpha-\ell+3j)(\alpha+\beta+\ell+j+2) \big)\displaybreak[0]\\
&+u^{i+j}\sum_{r}\textstyle
\binom{r}{i}\binom{r}{j-1}\binom{\beta+r-1}{r}\binom{\ell+k-1-r}{\ell-r}(1-u)^{\ell-r}3(\ell-j+1)(\beta+j-k),
\end{align*}
\begin{align*}
&\left(Z(Q_1-Q_0-2uQ_1)\right)_{ij}\displaybreak[0]\\ &
=u^{i+j}\sum_{r} \textstyle
\binom{r}{i}\binom{r}{j}\binom{\beta+r-1}{r}\binom{\ell+k-1-r}{\ell-r}(1-u)^{\ell-r}(-\alpha+\ell-3r)\displaybreak[0]\\
&\quad + u^{i+j}\sum_{r} \textstyle
\binom{r+1}{i}\binom{r+1}{j}\binom{\beta+r}{r+1}\binom{\ell+k-2-r}{\ell-r-1}(1-u)^{\ell-r}(3r+3j+2\alpha-2\ell+3),
\end{align*}
\begin{align*}
&\left(\textstyle{\frac{(\beta(1-u)-\alpha
u)}{u}}Z(Q_0+uQ_1)\right)_{ij}\\& \quad =u^{i+j}\sum_{r} \textstyle
\binom{r}{i}\binom{r}{j}\binom{\beta+r-1}{r}\binom{\ell+k-1-r}{\ell-r}
(1-u)^{\ell-r} \alpha(\ell-\alpha-3r)\displaybreak[0]
\\& \quad + u^{i+j}\sum_{r}\textstyle
\binom{r+1}{i}\binom{r+1}{j}\binom{\beta+r}{r+1}
\binom{\ell+k-2-r}{\ell-r+1}(1-u)^{\ell-r}(\alpha+\beta)(3r+\alpha-\ell+3),
\end{align*}
\begin{align*}
&((1-u)Z'(Q_1+uQ_0))_{i,j}\displaybreak[0]\\
&  =u^{i+j}\sum_{r} \textstyle
\binom{r}{i}\binom{r}{j}\binom{\beta+r-1}{r}\binom{\ell+k-2-r}{\ell-r+1}
(1-u)^{\ell-r}(r-\ell)(\alpha-\ell+3r)\displaybreak[0] \\
& +u^{i+j}\sum_{r}\textstyle
\binom{r+1}{i}\binom{r+1}{j}\binom{\beta+r}{r+1}\binom{\ell+k-2-r}{\ell-r+1}(1-u)^{\ell-r}\Big(3(r-j+1)(\ell-r+i+j)
\displaybreak[0]\\
&\quad \quad \quad \quad +(\alpha-\ell+3j)(\ell-r+i+j-1)\Big).
\end{align*}

By using the previous results we get that the identity \eqref{eqII}
is equivalent to
\begin{align*}
&\sum_{r=j}^{\ell}\!\!\textstyle
\binom{r}{i}\binom{r}{j}\binom{\beta+r-1}{r}\binom{\ell+k-1-r}{\ell-r}(1-u)^{\ell-r}
\textstyle{\frac{
3(r+1)(r+\beta-k+1)(\ell-r)(2r+2-i-j)}{(r-i+1)(r-j+1)}}\displaybreak[0]\\
&+(1-u)^{\ell-j+1}\textstyle\binom{j}{i}\binom{\beta+j-1}{j}
\binom{\ell+k-1-j}{\ell-j}3(\ell-j+k)(j-i) \displaybreak[0] \\
&-\!\!\sum_{r=j-1}^{\ell-1}\!\!\textstyle
\binom{r+1}{i}\binom{r+1}{j}\binom{\beta+r}{r+1}\binom{\ell+k-2-r}{\ell-r-1}(1-u)^{\ell-r}
3(\ell-r+k-1)(2r+2-i-j)\\ &=0,
\end{align*}
which easily follows.

In order to prove the identity \eqref{eqIII} we compute
\begin{align*}
&(ZV-V^*Z)_{ij}\\
&=(i-j)u^{i+j-1}\Big(-\sum_{r}\textstyle{\binom{r}{i}
\binom{r}{j}\binom{\beta+r-1}{r}}\binom{\ell+k-1-r}{\ell-r}
(1-u)^{\ell-r}(\alpha+\ell-r+1)\\
&+\sum_{r}\textstyle{
\binom{r}{i}\binom{r}{j}\binom{\beta+r-1}{r}\binom{\ell+k-1-r}{\ell-r}}(1-u)^{\ell-r+1}
(\alpha+\beta+i+j+\ell-r+1)\Big),
\end{align*}
\begin{align*}
&(P_1^*Z-ZP_1)_{ij}=(i-j)u^{i+j-1}\Big(\sum_{r}
\textstyle{\binom{r}{i}\binom{r}{j}\binom{\beta+r-1}{r}\binom{\ell+k-1-r}{\ell-r}}\\
& \qquad (1-u)^{\ell-r+1}(3i+4\alpha+3\beta+3k+6-3r+5\ell+3j)\displaybreak[0]\\
& -\sum_{r}
\textstyle{\binom{r}{i}\binom{r}{j}\binom{\beta+r-1}{r}\binom{\ell+k-1-r}{\ell-r}}
(1-u)^{\ell-r}(4\alpha+5\ell+3k+6-3r)\Big),
\end{align*}
\begin{align*}
&\Big(\big( \frac{\beta(1-u)-\alpha u}{u(1-u)}\big)((P_0^*+uP_1^*)Z-Z(P_0+uP_1))+ (P_0^*+uP_1^*)Z'\displaybreak[0]\\
&-Z'(P_0+uP_1)\Big)_{i,j} =u^{i+j-1}(i-j)\sum _{r}\textstyle{\binom
ri\binom rj \binom{\beta+r-1}{r}\binom{\ell+k-1-r}{\ell-r}}
(1-u)^{\ell-r}\displaybreak[0]\\
&\quad \quad\quad  \big(2( \alpha+2\ell+3k )(
\alpha+\ell-r+2)+2\alpha+\ell+6-3k-3r\big)\displaybreak[0]\\
&-u^{i+j-1}(i-j)\sum _{r}\textstyle\binom ri \binom r
j\binom{\beta+r-1}{r}
\binom{\ell+k-1-r}{\ell-r}(1-u)^{\ell-r+1}\big(\displaybreak[0]\\
&
2(\alpha+2\ell+3k)(\alpha+\beta+\ell-r+i+j+3)+3(i+j+\beta-3k+2-r-\ell)
\big ).
\end{align*}
Now it is easy to verify that \eqref{eqIII} is satisfied.

Finally the boundary conditions  \eqref{borde} can be easily check
and this concludes the proof of the theorem. \qed
\end{proof}

\section{The algebra of differential operators}\label{seccionalgebra}

Most of the results of this section are due to J. Tirao and they are
taken from  \cite{T2}.

Let $W=W(x)$ be a $L\times L$ matrix weight function  with finite
moments and let  $\{P_n\}$ be any sequence of matrix valued
orthogonal polynomials associated to a weight function $W$.

Let $$V_n=\{F\in M_{L\times L}(\CC)[x]:\deg(F)\leq n  \}$$  be the
set of all matrix valued polynomials in the variable $x$ of degree
less or equal to
  $n$.

\begin{prop}\label{decomposition}
    We have the following decomposition of $V_n$
  $$V_n=\bigoplus_{j=0}^n P_j^*M_{L\times L}(\CC) .$$
\end{prop}
\begin{proof} It is clear that $\sum_{j=0}^n P_j^* M_{L\times L}(\CC) $ is a
  subspace of $V_n$ and that for $n=0$ they are the same.
  Let us denote by $M_n$ the leading coefficient of $P_n^*$.
  If  $H=A_nx^n+A_{n-1}x^{n-1}+\cdots + A_0$ is a polynomial in
  $V_n$ then $H-P_n^*M_n^{-1}A_n$ is a
  polynomial of degree $\leq n-1$. Thus, by induction in $n$ we
  obtain that $H\in\sum_{j=0}^n P_j^*M_{L\times L}(\CC)$.
  In order to prove that this sum is a direct sum we assume that
  $P_0^*A_0^*+\cdots +P_n^*A_n^*=0$. By comparing, inductively the
  coefficients of $x^n, x^{n-1}, \dots x^0$ we obtain that $A_n=
  \cdots = A_0=0$. \qed
\end{proof}

\medskip
Let $\mathcal D$ be the algebra  of all differential operators of
the form
\begin{equation}\label{operadorordens}
  D= F_s(x) \frac{d^s}{dx^s}+F_{s-1}(x)
  \frac{d^{s-1}}{dx^{s-1}}+\cdots + F_1(x) \frac{d}{dx}+ F_0(x)
\end{equation}
whith  $F_j$ a polynomial function of degree less or equal to $j$.

\begin{thm}\label{Vn}
Let  $\{P_n\}$ be any sequence of matrix valued orthogonal
polynomials associated to $W$.
 If $D\in \mathcal D$ is symmetric  respect to
$W$ then $DP_n^*=P_n^*\Lambda_n$, for some matrix $\Lambda_n$.
\end{thm}
\begin{remark*}
We recall that $D$ is symmetric with respect to $W$ if $\langle DP,
Q\rangle=\langle P, DQ\rangle $, for all $P,Q$ polynomials. The
sequence $\{P_n\}$ is orthogonal with respect to $(\, , )$. The
bilinear forms $\langle \, , \rangle$ and $(\,,)$ are related  by
$\langle P,Q\rangle= (P^*,Q^*)^*$.
\end{remark*}
\begin{proof}
Since $D\in \mathcal D$  the operator $D$ preserves the vector
spaces $V_n$, for each $n\geq 0$.

For $n=0$ we have that $DP_0^*\in V_0$, thus
$DP_0^*=P_0^*\Lambda_0$. By induction we assume that
$DP_j^*=P_j^*\Lambda_j$, for each $0\leq j\leq n-1$. By Proposition
\ref{decomposition} we have that $DP_n^*=\sum_{i=0}^n P_i^*A_i$.
Thus, for each $0\leq j\leq n-1$ we have
$$\langle DP_n^*, P_j^*\rangle \textstyle
=\sum_{i=0}^j\langle P_i^*A_i, P_j^*\rangle= \sum_{i=1}^j
(P_i,P_j)^*A_i= (P_j,P_j)^*A_j.$$ On the other hand, since $D$ is
symmetric we obtain
$$\langle DP_n^*, P_j^*\rangle = \langle P_n^*, DP_j^*\rangle =
\langle P_n^*, P_j^*\Lambda_j\rangle = \left((P_n,
P_j)\Lambda_j\right)^* =0.$$ Thus $(P_j,P_j)^*A_j=0$  for each
$0\leq j\leq n-1$, which implies that $A_j=0$ because the matrix
$(P_j,P_j)$ is non singular. Therefore $DP_n^*=P_n^*\Lambda_n$ and
this concludes the proof. \qed
\end{proof}

\smallskip
Given   $\{P_n\}$ any sequence of matrix valued orthogonal
polynomials associated to the weight $W$, we define
\begin{equation}\label{algebraDW}
\mathcal D(W)=\{D\in \mathcal D: DP_n^*=P_n^* \Lambda_n(D),\forall
n\geq 0, \text{ for some  matrix } \Lambda_n(D)\}.
\end{equation}

\begin{prop}\label{propDW} We have
\begin{enumerate}
\item $\mathcal {D} (W)$ is a subalgebra of
$\mathcal D$ which does not depend on the sequence $\{P_n\}$.

\item  For each $n\in \NN_0$, the function $\Lambda_n:\mathcal
D(W)\longrightarrow M_{L\times L}( \CC)$
 given by $D\mapsto \Lambda_n(D)$
is a representation of the algebra $\mathcal D(W)$.

\item  The family $\{\Lambda_n\}_{n\geq 0}$ separates points of $\mathcal
D(W)$. That is, if $D_1$ and $D_2$ are distinct points of $\mathcal
D(W)$, then there exists $n_0\geq0$ such that
$\Lambda_{n_0}(D_1)\neq\Lambda_{n_0}(D_2)$.
\end{enumerate}
\end{prop}

\begin{proof}
It is easy to verify that $D(W)$ is a subalgebra of $\mathcal{D}$.
To prove that it is independent of the sequence $\{P_n\}$ we take
other sequence of orthogonal polynomials $\{Q_n\}$. Then
$Q_n=A_nP_n$, for some non singular matrix $A_n$. Then we have
$DQ_n^*=DP_n^*A_n^*= P_n^*\Lambda_n(D) A_n^*= Q_n^* \Upsilon_n(D),$
where $\Upsilon_n(D)=(A_n^*)^{-1}\Lambda_n(D) A_n^*$.

If $D_1$ and $D_2$ are in $D(W)$ then $$D_1D_2 P_n^*=D_1(P_n^*
\Lambda_n(D_2))=P_n^* \Lambda_n(D_1)\Lambda_n(D_2).$$ Hence
$\Lambda_n(D_1D_2)=\Lambda_n(D_1)\Lambda_n(D_2)$.

Let us assume that there exists $D\in D(W)$ such that
$\Lambda_n(D)=0$ for all $n\geq0$. To prove (3) we have to verify
that $D=0$. By hypothesis we have that $D=\sum_{i=0}^s F_i(x)
\frac{d^i}{dx^i}$
 satisfies  $DP_n^*=0$, for all $n\geq 0$. For $n=0$ we
obtain $F_0P_0^*=0$, thus $F_0=0$. \\
By induction, we may assume that $F_i=0$ for $0\leq i\leq j-1$, with
$j\leq s$. Then $0=DP_j^*=\sum_{i=1}^j
F_i(x)\frac{d^i(P_j^*)}{dx^i}=F_j(x)j!M_j$, where $M_j$ is the
leading coefficient of $P_j$, which is non singular. Therefore
$F_j=0$. This concludes the proof.  \qed
\end{proof}

\begin{cor}\label{autovconmutan} The operators $D_1$ and $D_2$ in the algebra $\mathcal
D(W)$ commute if and only if the matrices $\Lambda_n(D_1)$
 and $\Lambda_n(D_2)$ commute for all $n\in \NN_0$.
\end{cor}
\begin{proof}
  By Proposition \ref{propDW}, (3) we have that $D_1D_2=D_2D_1$ if and only if
  $\Lambda_n(D_1D_2)=\Lambda_n(D_2D_1)$ for all $n$. From Proposition \ref{propDW}, (2) we get
    $\Lambda_n(D_1D_2)=\Lambda_n(D_1)\Lambda(D_2)$. \qed
\end{proof}

\smallskip

\begin{prop}\label{Pro2} Let $\{Q_n\}$ be the sequence of
monic orthogonal polynomials.  Let $D=\sum_{i=0}^s
F_i(u)\frac{d^i}{du^i}$ such that $DQ_n^*= Q_n^*\Gamma_n.$ Then
\begin{equation}\label{lambda}
\Gamma_n=\sum_{0\leq i\leq s }[n]_i A_i^i\quad\quad \text{for
all}\quad n\ge0,
\end{equation}
where $A_i^i$ is the coefficient of $x^i$ in the polynomial $F_i$.
\end{prop}
\begin{remark*}
Here we are using the notation $[n]_i=n(n-1)\dots (n-i+1)$ for
$n\geq 1$, and $[n]_0=1$, for $n\geq 0$.
\end{remark*}

\begin{proof}
From
$$ \sum_{0\leq i\leq s}\textstyle F_i(u)\frac{d^i Q_n^*}{du^i}(u)=Q_n^*(u)\Gamma_n ,$$
by comparing the monomials of degree $n$ we get $\sum_{0\leq i\leq s
} [n]_i A_i^i=\Gamma_n.$ \qed
\end{proof}

\begin{remark}
 Observe that in particular, Proposition \ref{Pro2} implies that the eigenvalues
 $\Gamma_n$ is a polynomial function on $n$ of degree less or equal
 to $\deg(D)$.
\end{remark}
\section{The operator $E$}

\subsection{$D$ and $E$ commute}\label{seccDEconmutan} In this
subsection we use the results described in Section
\ref{seccionalgebra} to give an elegant proof of the fact that the
operators $D$ and $E$ commute. Of course we also can verify this by
making the explicit computations.

\begin{thm}\label{comnutan}
The differential operators $D$ and $E$, introduced respectively in
Theorems \ref{D1salto} and \ref{E1salto}, commute.
\end{thm}

\begin{proof}
From Theorem \ref{E1salto} the operator $E$ is symmetric with
respect to the weight $W$. Thus $E$ belongs to the algebra
$\mathcal{D}(W)$ defined in \eqref{algebraDW} (See Theorem
\ref{Vn}). To see that $D$ and $E$ commute it is enough to verify
that the corresponding eigenvalues commute. (See Corollary
\ref{autovconmutan}).

Let $\{Q_n\}$ be the monic sequence of orthogonal polynomials. Then
for any $D\in \mathcal{D}(W)$, we have $DQ_n^*=Q_n^*\Gamma_n(D)$,
where the eigenvalue $\Gamma_n(D)$ is given  explicitly in terms of
the coefficients of the differential operator $D$ (see Proposition
\ref{Pro2}).

For the operators $D$ and $E$ introduced respectively in Theorems
\ref{D1salto} and \ref{E1salto}, these eigenvalues are
\begin{align*}
  \Gamma_n(D)&= - n(U+n-1)-V \\
  \Gamma_n(E)&= -n(n-1)Q_1+nP_1-(\alpha+2\ell+3k)V,
\end{align*}
where the matrices $U,V,Q_1$ and $P_1$ are given in Theorems
\ref{D1salto} and \ref{E1salto}. Explicitly we have
\begin{align*}
  \Gamma_n(D)&= -\textstyle{\sum_{i=0}^\ell
\bigl(
n(n+\alpha+\beta+\ell+i+1)+i(i+\alpha+\beta-k+1)\bigr)E_{ii}}\\
&\quad +\textstyle{\sum_{i=0}^{\ell-1}
(\ell-i)(\beta+i-k+1)E_{i,i+1}}\displaybreak[0]\\
\Gamma_n(E)&= -\textstyle{\sum_{i=0}^\ell \bigl(
n(\alpha-\ell+3i)(n+\alpha+\beta+\ell+i+1)}\\
&\textstyle{\qquad \qquad\qquad +(\alpha+2\ell+3k)i(i+\alpha+\beta-k+1)\bigr)E_{ii}}\\
&\quad +\textstyle{\sum_{i=0}^{\ell-1} (\ell-i)(\beta+i-k+1)(\alpha
+ 2\ell+3k+3n)E_{i,i+1}}.
\end{align*}

Now it is easy to verify that
\begin{equation}\label{conmutan1}
\Gamma_n(E)= (\alpha+2\ell+3k+3n) \Gamma_n(D)+
3n(\ell+k+n)(n+\alpha+\beta+\ell+1)I. \end{equation}
Thus the matrix
$\Gamma_n(E)$ commutes with  $\Gamma_n(D)$ and by Corollary
\ref{autovconmutan} we have that $D$ and $E$ commute. \qed
\end{proof}

\subsection{The eigenfunctions of $E$.}

In Subsection \ref{orthogsubseccion} we give a sequence $\{P_w\}_w$
of matrix valued polynomials, which are orthogonal  with respect to
$W$ and  eigenfunctions of the differential operator $D$. The rows
$P_w^j$ of $P_w$ are orthogonal polynomials of degree $w$ and they
satisfy $DP_w^j=\lambda_j(w) P_w^j$.

 Since $D$ and $E$ commute, it follows
that $E$ preserves the eigenspaces of $D$. Therefore if an
eigenvalue $\lambda=\lambda_j(w)$ has multiplicity one, then the
vector valued polynomial $P_w^j$ is also an eigenfunction of the
differential operator $E$. In the next theorem, we prove that this
is true, even if the multiplicity of an eigenvalue is bigger than
one.

\begin{thm}\label{orthopolyE}
  The sequence $\{P_w\}_w$ of orthogonal polynomials associated to
  the  pair $\{W,D\}$ satisfies
  $$EP_w^*(u)=P_w^*(u) \Lambda_w(E),$$
  where $\Lambda_w(E)=\displaystyle\sum_{0\leq j\leq \ell} \mu_j(w) E_{jj}$,
  and
  $$\mu_j(w)=-w(w+\alpha+\beta+\ell+j+1)(\alpha-\ell+3j)-j(j+\alpha+\beta-k+1)(\alpha+2\ell+3k).$$
\end{thm}

\begin{proof}
Let $\{Q_w^*\}_{w\geq0}$ be the sequence of monic orthogonal
polynomials. Since $E$ is symmetric with respect to the weight $W$,
Theorem \ref{Vn} says that $EQ_w^*=Q_w^*\Gamma_w(E)$ for some matrix
$\Gamma_w(E)$. If $Q_w^*=P_w^*A_w^*$ then we have that
$DP_w^*A_w^*=P_w^*A_w^*\Gamma_w(D)$ and
$EP_w^*A_w^*=P_w^*A_w^*\Gamma_w(E)$. Therefore
\begin{align*}
\Lambda_w(D)&=A_w^*\Gamma_w(D)(A_w^*)^{-1},\\
\Lambda_w(E)&=A_w^*\Gamma_w(E)(A_w^*)^{-1}.
\end{align*}
Thus from \eqref{conmutan1} we obtain that
$$\Lambda_w(E)=(\alpha+2\ell+3k+3w)
\Lambda_w(D)+ 3w(\ell+k+w)(w+\alpha+\beta+\ell+1)I.$$

Observe that the fact that $\Lambda_w(D)$ is a diagonal matrix
implies that $\Lambda_w(E)$ is diagonal. Moreover the eigenvalue
$\mu_j(w)=(\Lambda_w(E))_{jj}$ is given by
\begin{align*}
\mu_j(w)&=(\alpha+2\ell+3k+3w)(-w(w+\alpha+\beta+\ell+j+1)\\
&\quad -j(j+\alpha+\beta-k+1))+3w(\ell+k+w)(w+\alpha+\beta+\ell+1)\\
&=-w(w+\alpha+\beta+\ell+j+1)(\alpha-\ell+3j)\\
&\quad-j(j+\alpha+\beta-k+1)(\alpha+2\ell+3k).
\end{align*}
This concludes the proof of the theorem. \qed
\end{proof}

\subsection{The operator algebra generated by $D$ and $E$}\label{algebra}
In this subsection we study the algebra generated by the
differential operators $D$ and $E$.

Let $\CC[x,y]$ be the algebra of all polynomials in the variables
$x$ and $y$ with complex coefficients.

\begin{thm}\label{alggenDE}
The algebra of differential operators generated by $D$ and $E$  is
isomorphic to the  quotient algebra $\CC[x,y]/ \langle Q \rangle$,
where $\langle Q \rangle$ denotes the ideal generated by the
polynomial
$$Q(x,y)=\prod_{j=0}^\ell \left(y-(\alpha-\ell+3j)x+3j(\ell-j+k)(j+\alpha+\beta-k+1) \right).$$
\end{thm}
\begin{proof}
The algebra of differential operators generated by $D$ and $E$ is
isomorphic to the quotient algebra $\CC[x,y]/I$ where $I=\{
p\in\CC[x,y]:p(D,E)=0\}$.

Since $\Lambda_w$ is a representation which separates points of
$\mathcal{D}(W)$ (Proposition \ref{propDW}), we have that $p(D,E)=0$
if and only if
$$\Lambda_w(p(D,E))=p(\Lambda_w(D),\Lambda_w(E))=0,\text{ for all }w.$$
Moreover, since the matrices $\Lambda_w(D)$ and $\Lambda_w(E)$ are
diagonal matrices, we have that $p(\Lambda_w(D),\Lambda_w(E))=0$ if
and only if $p((\Lambda_w(D))_{jj},(\Lambda_w(E))_{jj})=0$ for all
$0\leq j\leq\ell$. Thus the ideal $I$ is
$$I=\{p\in\CC[x,y]:p(\lambda_j(w),\mu_j(w))=0, \text{for } j=0,1,\dots  \ell\}.$$

 Let $p_j(x,y)$ be the polynomial
$$p_j(x,y)=y-(\alpha-\ell+3j)x+3j(\ell-j+k)(j+\alpha+\beta-k+1).$$
It is easy to verify that $p_j(\lambda_j(w),\mu_j(w))=0,$ for all
$w\geq0$. Therefore $Q(x,y)=\prod_{j=0}^\ell p_j(x,y)$ belongs to
the ideal $I$.

 On the other hand we have that any $f\in I$  vanishes in all
points of the form $(x,y)$ with
$y=(\alpha-\ell+3j)x+3j(\ell-j+k)(j+\alpha+\beta-k+1)$, for each
$j=0,\dots ,\ell$. In fact if we let,
$$a_j=\alpha-\ell+3j \qquad  b_j=-3j(\ell-j+k)(j+\alpha+\beta-k+1)\quad (j=0,1,\dots ,\ell)$$
then we observe  that
 the polynomial $f(x,a_jx+b_j)$ has infinitely many roots,  because   $f(\lambda_j(w),\mu_j(w))=0$
and $\mu_j(w)=a_j\lambda_j(w)+b_j$.\\
 Any polynomial in $\CC[x,y]$ is also a polynomial in $x$ and $y-ax-b$. Then it is clear that if
$p(x,y)=0$ in the line $y=ax+b$ then $p$ is divisible by $y-ax-b$.
\\
Thus we have that if $f$ belongs to the ideal $I$ then $f\in
\cap_{j=0}^\ell \langle p_j \rangle= \langle \textstyle \prod_j
p_j\rangle.$ Therefore we have that the ideal $I$ is generated by
the polynomial $Q(x,y)$, which concludes the proof of the
Theorem.\qed
\end{proof}

\end{document}